# Abel-Grassmann's Groupoids of Modulo Matrices

AMANULLAH[1], M. RASHAD[1], I. AHMAD[1]

## ABSTRACT

The binary operation of usual addition is associative in all matrices over $R$. However, a binary operation of addition in matrices over $Z_n$ of a nonassociative structures of AG-groupoids and AG-groups are defined and investigated here. It is shown that both these structures exist for every integer $n \geq 3$. Various properties of these structures are explored like: (i). Every AG-groupoid of matrices over $Z_n$ is transitively commutative AG-groupoid and is a cancellative AG-groupoid if n is prime. (ii). Every AG-groupoid of matrices over $Z_n$ of type-II is a $T^3$-AG-groupoid. (iii). An AG-groupoid of matrices over $Z_n; G_{nAG}(t,u)$ is an AG-band, if $t+u=1(\bmod n)$.

Keywords: AG-groupoid and AG-group of matrices over $Z_n$; $T^3$-AG-groupoid; Transitively commutative AG-groupoid; Cancellative AG-groupoid.

## 1. INTRODUCTION

A magma that satisfies the left invertive law: $(ab)c=(cb)a$ is called left almost semigroup abbreviated as LA-semigroup [1], or Abel-Grassmann's groupoid abbreviated as AG-groupoid [2]. An AG-group $G$ is an AG-groupoid which has the left identity and has inverse of each of its elements. Both these structures are nonassociative in general, and so one has to play the game of brackets in a defined way. In an AG-groupoid $G$ an element $a \in G$ is called idempotent if $a=a^2$, and $G$ is called idempotent or AG-2-band or simply AG-band if each of its elements is idempotent [3]. AG-groupoids generalize commutative semigroups, while an AG-group generalize an Abelian group. These structures have a variety of applications in geometry, flocks theory, topology, finite mathematics and many more [4-7]. Many new classes of AG-groupoids have recently been introduced and characterized [8-12]. Fuzzification of AG-groupoids and AG-groups has also been done see for instance [13,14]. M. Shah studied a lot about the AG-groupoids and AG-groups exclusively [4]. However, the construction of these structures, as well as of other algebraic structures remains a difficult job for the researchers. In this article we try to introduce construction of these structures especially in matrices. A matrix $A$ over a field $F$ is a rectangular array of scalars represented by:

$$A=\begin{bmatrix} a_{11} & a_{12} & \cdots & a_{1q} \\ a_{21} & a_{22} & \cdots & a_{2q} \\ \vdots & \vdots & \vdots & \vdots \\ a_{p1} & a_{p2} & \cdots & a_{pq} \end{bmatrix}$$

The rows of matrix $A$ are the $p$ horizontal entries $(a_{11},a_{12},\cdots,a_{1q})$, $(a_{21},a_{22},\cdots,a_{2q})$, $\cdots,(a_{p1},a_{p2},\cdots,a_{pq})$, and the column of $A$ are the $q$ vertical entries

$$\begin{bmatrix} a_{11} \\ a_{21} \\ \vdots \\ a_{p1} \end{bmatrix}, \begin{bmatrix} a_{12} \\ a_{22} \\ \vdots \\ a_{p2} \end{bmatrix} \text{ and } \begin{bmatrix} a_{1q} \\ a_{2q} \\ \vdots \\ a_{pq} \end{bmatrix}$$

Note that $a_{ij}$ represents the entry in $i^{th}$ row and $j^{th}$ column. A matrix with $p$ rows and $q$ columns is called $p \times q$ matrix. A matrix having only one row is called row matrix or row vector, and a matrix having only one column is called a column matrix or column vector.

Let $A=\left[a_{ij}\right]$ and $B=\left[b_{ij}\right]$ be two matrices of the same order, say $p\times q$. Then $A+B$, is a matrix obtained by adding the corresponding entries from $A$ and $B$ i.e. $A+B=\left[a_{ij}+b_{ij}\right]$. The product of the matrix $A$ by a scalar $c$, $cA$ is a matrix obtained by multiplying each entry of $A$ by $c$ i.e. $cA=\left[ca_{ij}\right]$.

W. B. Vasantha [15] has defined some nonassociative construction for groupoids on matrices as follows:

**Definition-1.** Let $G_1=\{(x_1,x_2,\cdots,x_q)\mid x_i\in Z_n;1\le i\le q\};n\ge 3$ be the collection of row matrices over $Z_n$. Define a binary operation '*' on $G_1$ as follows:

$$\begin{aligned}\left(x_1,x_2,\cdots,x_q\right)*\left(y_1,y_2,\cdots,y_q\right) &\cong \left[t\left(x_1,x_2,\cdots,x_q\right)+u\left(y_1,y_2,\cdots,y_q\right)\right]\ (\text{mod } n)\\ &\cong \left[\left(tx_1+uy_1\right)(\text{mod } n),\ \cdots,\ \left(tx_q+uy_q\right)(\text{mod } n)\right]\end{aligned}$$

where $t,\ u\in Z_n-\{0\}$, $t=u$ and $(t,\ u)=1$ for all $\left(x_1,\ x_2,\ \cdots,\ x_q\right)$, $\left(y_1,\ y_2,\ \cdots,\ y_q\right)\in G_1$. Thus $(G_1(t,\ u),*)$ is a groupoid of row matrix over $Z_n$.

**Definition-2.** Let $G_2=\{(x_1,\ x_2,\ \cdots,\ x_p)^t\mid x_i\in Z_n;\ 1\le i\le p\};\ n\ge 3$ be the collection of $p\times 1$ column matrices over $Z_n$. Choose $t,\ u\in Z_n-\{0\}$, $t=u$ and $(t,\ u)=1$ for all $\left(x_1,\ x_2,\ \cdots,\ x_p\right)^t$, $\left(y_1,\ y_2,\ \cdots,\ y_p\right)^t\in G_2$. Define '*' on $G_2$ as follows:

$$\left(x_1,\ x_2,\ \cdots,\ x_p\right)^t*\left(y_1,\ y_2,\ \cdots,\ y_p\right)^t=\begin{bmatrix}x_1\\x_2\\\vdots\\x_p\end{bmatrix}*\begin{bmatrix}y_1\\y_2\\\vdots\\y_p\end{bmatrix}\cong\begin{bmatrix}(tx_1+uy_1)(\text{mod } n)\\(tx_2+uy_2)(\text{mod } n)\\\vdots\\(tx_p+uy_p)(\text{mod } n)\end{bmatrix}$$

Thus $(G_2(t,\ u),*)$ is a groupoid of column matrix over $Z_n$.

**Definition-3.** Let $G=\{[m_{ij}]\mid 1\le i\le p,\ 1\le j\le q\}$, be the collection of $p\times q$ matrices over $Z_n;\ n\ge 3$. For non zero distinct integers $t,\ u$ in $Z_n$, define '*' on $G$, for two matrices $M=[m_{ij}]$ and $N=[n_{ij}]$ in $G$ as follow;

$$M*N=\left[m_{ij}\right]*\left[n_{ij}\right]\cong\left[\left(tm_{ij}+un_{ij}\right)(\text{mod } n)\right]$$

Thus $(G(t,\ u),*)$ is any groupoid on $p\times q$ matrix over $Z_n$. Moreover, we generally represent this groupoid of (row matrix, column matrix or $p\times q$ matrix) over $Z_n$ by $G_n(t,\ u)$. We will use $G(n)$ to denote the class of $G_n(t,\ u)$ for distinct integers $t,\ u\in Z_n-\{0\}$, and $(t,\ u)=1$, that is:

$$G(n)=\{G_n(t,\ u),\text{ for distinct integers } t,u\in Z_n-\{0\},\text{ and } (t,\ u)=1\}.$$

By putting some additional conditions on $t$ and $u$ in $G(n)$, we get some other classes of groupoids of matrices over $Z_n$ as follow:

(i) Type-I, if for non-zero distinct integers $t,u$ in $Z_n$ and $(t,\ u)=1$.

(ii) Type-II, if $t,\ u \in Z_n - \{0\}$, such that $t = u$.

(iii) Type-III, if $t,\ u \in Z_n$, where $t$ or $u$ is zero.

In the following section we show the existence of AG-groupoid of matrices over $Z_n$, and find its relations with some of the already known classes of AG-groupoids.

## 2. EXISTENCE OF AG-GROUPOID OF MATRICES OVER $Z_n$

The following theorem shows the existence AG-groupoid of matrices over $Z_n$ where $n \geq 3$, and indeed it introduces a simple way of construction of these AG-groupoids of any finite order.

**Theorem-1.** $G_n(t,\ u)$, is an AG-groupoid of matrices over $Z_n$, if $t^2 \cong u \pmod{n}$ for any $t,\ u \in Z_n$.

**Proof.** Let $G_n(t,\ u)$, satisfies $t^2 \cong u \pmod{n}$ for any $t,\ u \in Z_n$. To show that $G_n(t,\ u)$ is an AG-groupoid of row matrices (column matrices or $p \times q$ matrices), it is sufficient if we show the left invertive law; $(A*B)*C = (C*B)*A; \forall A,\ B,\ C \in G_n(t,\ u)$ holds;

$$(A*B)*C = \left(\left[a_{ij}\right]*\left[b_{ij}\right]\right)*\left[c_{ij}\right] \cong \left[\left(ta_{ij} + ub_{ij}\right)\pmod{n}\right]*\left[c_{ij}\right]$$

$$\Rightarrow (A*B)*C \cong \left(t^2 a_{ij} + tub_{ij} + uc_{ij}\right)\pmod{n} \tag{1}$$

and

$$(C*B)*A = \left(\left[c_{ij}\right]*\left[b_{ij}\right]\right)*\left[a_{ij}\right] \cong \left(t^2 c_{ij} + tub_{ij} + ua_{ij}\right)\pmod{n} \tag{2}$$

also

$$A*(B*C) = \left[a_{ij}\right]*\left(\left[b_{ij}\right]*\left[c_{ij}\right]\right) \cong \left(ta_{ij} + utb_{ij} + u^2 c_{ij}\right)\pmod{n} \tag{3}$$

This implies that $G_n(t,\ u)$ is nonassociative AG-groupoid of matrices over $Z_n$ by (1-3).

We denote this AG-groupoid of matrices over $Z_n$ by $G_{nAG}(t,\ u)$, and $G_{AG}(n)$ will represent the class that contains all AG-groupoids of matrices over $Z_n$. Now by varying values of $t$ and $u$ and by imposing some additional conditions on $t$ and $u$, we get different classes of AG-groupoid of matrices over $Z_n$ for some fixed integer $n \geq 3$. Thus the so obtained new classes of AG-groupoids of matrices over $Z_n$ will be denoted by $G_{AG-I}(n)$, $G_{AG-II}(n)$ and $G_{AG-III}(n)$. The following example shows the existence of these AG-groupoids of matrices over $Z_n$.

**Example-1.** $G_3(2,\ 1)$ is an AG-groupoid of matrices over $Z_3$, that is $G_{3AG}(2,\ 1) \in G_{AG}(3)$.

**Solution.** As $G_3(2,\ 1) \in G(3)$, to show that it is an AG-groupoid of matrices over $Z_3$, that is, $G_{3AG}(2,\ 1) \in G_{AG}(3)$, we show that it satisfies the left invertive law:

$$(A*B)*C = \left(\left[a_{ij}\right]*\left[b_{ij}\right]\right)*\left[c_{ij}\right] \cong \left[\left(2a_{ij} + b_{ij}\right)\pmod{n}\right]*\left[c_{ij}\right]$$

$$\Rightarrow (A*B)*C \cong \left[\left(a_{ij} + 2b_{ij} + c_{ij}\right)\pmod{n}\right] \tag{4}$$

and

$$(C*B)*A = \left(\left[c_{ij}\right]*\left[b_{ij}\right]\right)*\left[a_{ij}\right] \cong \left[\left(c_{ij} + 2b_{ij} + a_{ij}\right)\pmod{n}\right] \tag{5}$$

also

$$A*(B*C) = \left[a_{ij}\right]*\left(\left[b_{ij}\right]*\left[c_{ij}\right]\right) \cong \left[\left(2a_{ij} + 2b_{ij} + c_{ij}\right)(\text{mod } n)\right] \quad (6)$$

This implies that $G_{3AG}(2,\ 1) \in G_{AG}(3)$ by (4, 5), and is nonassociative by (4,6).

**Example-2.** $G_8(6,\ 4)$ is an AG-groupoid of matrices over $Z_8$ of Type-I, that is

$$G_{8AG}(6,\ 4) \in G_{AG-I}(8).$$

The following examples show some various types of AG-groupoids of matrices over $Z_n$ for $n \geq 3$.

**Example-3.** $G_{AG}(3) = \{G_{3AG}(2,1)\}$ and $G_{AG-II}(3) = \{G_{3AG}(1,1)\}$.

**Example-4.** $G_{AG}(4) = \{G_{4AG}(3,1)\}$, $G_{AG-II}(4) = \{G_{4AG}(1,1)\}$ and $G_{AG-III}(4) = \{G_{4AG}(2,0)\}$.

**Example-5.** $G_{AG}(5) = \{G_{5AG}(4,1),\ G_{5AG}(3,4)\}$, $G_{AG-I}(5) = \{G_{5AG}(2,4)\}$ and $G_{AG-II}(5) = \{G_{5AG}(1,1)\}$.

**Example-6.** $G_{AG}(6) = \{G_{6AG}(5,1)\}$, $G_{AG-I}(6) = \{G_{6AG}(2,4)\}$ and $G_{AG-II}(6) = \{G_{6AG}(1,1),\ G_{6AG}(3,3),\ G_{6AG}(4,4)\}$ and so on.

The following corollaries are straight away by Theorem-1.

**Corollary-1.** Any $G_{nAG}(t,\ u) \in G_{AG-II}(n)$ is an abelian group (an AG-group), if $t = u = 1$.

**Corollary-2.** Any $G_{nAG}(t,\ u) \in G_{AG-II}(n)$ is a commutative semigroup, if $t = u \neq 1$.

**Example-7.** AG-groupoids of matrices $G_{6AG}(3,\ 3)$ and $G_{6AG}(4,\ 4)$ over $Z_6$ are commutative semigroups in $G_{AG-II}(6)$.

## 3. SOME PROPERTIES OF AN AG-GROUPOID OF MATRICES OVER $Z_n$

In this section we investigate the relations of these new classes of AG-groupoid of matrices over $Z_n$ with some of the already known classes of AG-groupoids that includes the following:

(i). $T^3$-AG-groupoid, if it is;

  (a). $T_l^3$-AG-groupoid, that is, if $a*b = a*c \Rightarrow b*a = c*a$ [16].

  (b). $T_r^3$-AG-groupoid, that is, if $b*a = c*a \Rightarrow a*b = a*c$ [16].

(ii). Transitively commutative AG-groupoid, if $a*b = b*a$ and $b*c = c*b$ $\Rightarrow a*c = c*a$ [16].

(iii). Cancellative AG-groupoid, an element $a$ in an AG-groupoid $G$ is left (right) cancellative, if $a \cdot x = a \cdot y \Rightarrow x = y$ $(x \cdot a = y \cdot a \Rightarrow x = y)$, and $G$ is left (right) cancellative if each of its elements is left (right) cancellative [17]

(iv). AG-band, if $a*a = a$ [3].

**Theorem-2.** Every AG-groupoid of matrices over $Z_n$ of type-II that is $G_{AG-II}(n)$ is a $T^3$-AG-groupoid.

**Proof.** To show that $G_{AG-II}(n)$ is a $T^3$-AG-groupoid, it is sufficient if we show that an arbitrary AG-groupoid of matrices over $Z_n$ of type-II is $T_l^3$-AG-groupoid and $T_r^3$-AG-groupoid.

For $T_l^3$-AG-groupoid, let $A,\ B,\ C \in G_{AG-II}(n)$, and

$$A * B = A * C$$
$$\Rightarrow \left[a_{ij}\right] * \left[b_{ij}\right] = \left[a_{ij}\right] * \left[c_{ij}\right]$$
$$\Rightarrow \left[\left(ta_{ij} + ub_{ij}\right)\left(\text{mod } n\right)\right] \cong \left[\left(ta_{ij} + uc_{ij}\right)\left(\text{mod } n\right)\right]$$
$$\Rightarrow \left[\left(ub_{ij}\right)\left(\text{mod } n\right)\right] \cong \left[\left(uc_{ij}\right)\left(\text{mod } n\right)\right]$$
$$\left[\left(tb_{ij}\right)\left(\text{mod } n\right)\right] \cong \left[\left(tc_{ij}\right)\left(\text{mod } n\right)\right] \quad (\because t = u) \tag{7}$$

Now

$$B * A = \left[b_{ij}\right] * \left[a_{ij}\right]$$
$$\cong \left[\left(tb_{ij} + ua_{ij}\right)\left(\text{mod } n\right)\right]$$
$$\cong \left[\left(tc_{ij} + ua_{ij}\right)\left(\text{mod } n\right)\right] \text{ by } (7)$$
$$= \left[c_{ij}\right] * \left[a_{ij}\right]$$
$$\Rightarrow B * A = C * A.$$

Hence $G_{AG-II}\ (n)$ is $T_l^3$-AG-groupoid. Similarly we can show that $G_{AG-II}\ (n)$ is $T_r^3$-AG-groupoid. Hence $G_{AG-II}\ (n)$ is $T^3$-AG-groupoid.

**Theorem-3.** Every AG-groupoid of matrices over $Z_p$ that is $G_{pAG}\,(t,\ u)$ is a $T^3$-AG-groupoid, if $p$ is prime and $u \in Z_p - \{0\}$.

**Proof.** To show that every $G_{pAG}\,(t,\ u)$ is a $T^3$-AG-groupoid for $u \in Z_p - \{0\}$ it issufficient to show that, $G_{pAG}\,(t,\ u)$ is $T_l^3$-AG-groupoid and $T_r^3$-AG-groupoid.

For $T_l^3$-AG-groupoid, let $A, B, C \in G_{pAG}\,(t,\ u)$, and

$$A * B = A * C$$
$$\Rightarrow \left[a_{ij}\right] * \left[b_{ij}\right] = \left[a_{ij}\right] * \left[c_{ij}\right]$$
$$\Rightarrow \left[\left(ta_{ij} + ub_{ij}\right)\left(\text{mod } p\right)\right] \cong \left[\left(ta_{ij} + uc_{ij}\right)\left(\text{mod } p\right)\right]$$
$$\Rightarrow \left[u\left(b_{ij} - c_{ij}\right)\left(\text{mod } p\right)\right] \cong \left[0\left(\text{mod } p\right)\right],$$

as $p \nmid u$, because a non-zero $u$ and $\left(b_{ij} - c_{ij}\right)$ both are less then $p$, where $p$ is prime. Therefore, $p \mid \left[b_{ij} - c_{ij}\right] \Rightarrow \left[b_{ij}\,(\text{mod } p)\right] \cong \left[c_{ij}\,(\text{mod } p)\right]$, consequently;

$$B * A = \left[b_{ij}\right] * \left[a_{ij}\right]$$
$$\cong \left[\left(tb_{ij} + ua_{ij}\right)\left(\text{mod } p\right)\right]$$
$$\cong \left[\left(tc_{ij} + ua_{ij}\right)\left(\text{mod } p\right)\right]$$
$$= \left[c_{ij}\right] * \left[a_{ij}\right]$$
$$\Rightarrow B * A = C * A.$$

Hence every AG-groupoid of matrices over $Z_p$ is a $T_l^3$-AG-groupoid. Similarly, we can show that every AG-groupoid of matrices over $Z_p$ is a $T_r^3$-AG-groupoid. Hence for non zero $u \in Z_p$, every AG-groupoid of matrices over $Z_p$ is a $T^3$-AG-groupoid.

**Example-8.** $G_{5AG}(3,\ 4)$ is a $T^3$-AG-groupoid. However, the result is not true in general. For example; $G_{8AG}\ (6,\ 4)$ is not a $T^3$-AG-groupoid.

From the following theorem it is clear that the collection of AG-groupoids of matrices over $Z_n$ in any class is a subclass of transitively commutative AG-groupoid.

**Theorem-4.** Every $G_{nAG}(t,\ u)$ is a transitively commutative AG-groupoid.

**Proof.** To show that every $G_{nAG}(t,\ u)$ is a transitively commutative AG-groupoid, we show that an arbitrary AG-groupoid of matrices is transitively commutative AG-groupoid. Let $A,\ B,\ C \in G_{nAG}(t,\ u)$, and

$$A*B = B*A$$
$$\Rightarrow [a_{ij}]*[b_{ij}] = [b_{ij}]*[a_{ij}]$$
$$\Rightarrow [(ta_{ij} + ub_{ij})(\text{mod } n)] \cong [(tb_{ij} + ua_{ij})(\text{mod } n)]$$
$$\Rightarrow [(t(a_{ij} - b_{ij}) + u(b_{ij} - a_{ij}))(\text{mod } n)] \cong [0_{ij}\ (\text{mod } n)] \tag{9}$$

also,

$$B*C = C*B$$
$$\Rightarrow [b_{ij}]*[c_{ij}] = [c_{ij}]*[b_{ij}]$$
$$\Rightarrow [(tb_{ij} + uc_{ij})(\text{mod } n)] \cong [(tc_{ij} + ub_{ij})(\text{mod } n)]$$
$$\Rightarrow [(t(b_{ij} - c_{ij}) + u(c_{ij} - b_{ij}))(\text{mod } n)] \cong [0_{ij}\ (\text{mod } n)], \tag{10}$$

as $n \,|\, [t(a_{ij} - b_{ij}) + u(b_{ij} - a_{ij})]$ and $n \,|\, [t(b_{ij} - c_{ij}) + u(c_{ij} - b_{ij})]$ by (9,10)

$$\Rightarrow n \mid [t(a_{ij} - b_{ij}) + u(b_{ij} - a_{ij}) + t(b_{ij} - c_{ij}) + u(c_{ij} - b_{ij})]$$
$$\Rightarrow n \mid [t(a_{ij} - c_{ij}) + u(c_{ij} - a_{ij})]$$
$$\Rightarrow [(t(a_{ij} - c_{ij}) + u(c_{ij} - a_{ij}))(\text{mod } n)] \cong [0_{ij}\ (\text{mod } n)]$$
$$\Rightarrow [((ta_{ij} + uc_{ij}) - (tc_{ij} + ua_{ij}))(\text{mod } n)] \cong [0_{ij}\ (\text{mod } n)]$$
$$\Rightarrow [(ta_{ij} + uc_{ij})(\text{mod } n)] \cong [(tc_{ij} + ua_{ij})(\text{mod } n)]$$
$$\Rightarrow [a_{ij}]*[c_{ij}] = [c_{ij}]*[a_{ij}]$$
$$\Rightarrow A*C = C*A.$$

Hence every $G_{nAG}(t,\ u)$ is a transitively commutative AG-groupoid.

**Example-9.** $G_{7AG}\ (5,\ 4)$ is transitively commutative AG-groupoid.

**Theorem-5.** Every $G_{pAG}(t,\ u)$ is a cancellative AG-groupoid, if $p$ is prime and $u \in Z_p - \{0\}$.

**Proof.** Let $p$ is prime and $u \in Z_p - \{0\}$, to show that $G_{pAG}(t, u)$ is a cancellative AG-groupoid, we show that it is left cancellative AG-groupoid and right cancellative AG-groupoid.
For left cancellativity, let

$$\begin{aligned} A * X &= A * Y \\ \Rightarrow \left[a_{ij}\right] * \left[x_{ij}\right] &= \left[a_{ij}\right] * \left[y_{ij}\right] \\ \Rightarrow \left[\left(ta_{ij} + ux_{ij}\right)(\text{mod } p)\right] &\cong \left[\left(ta_{ij} + uy_{ij}\right)(\text{mod } p)\right] \\ \Rightarrow \left[u\left(x_{ij} - y_{ij}\right)(\text{mod } p)\right] &\cong \left[0_{ij}\,(\text{mod } p)\right] \end{aligned}$$

$\Rightarrow p \mid u\left(x_{ij} - y_{ij}\right)$ as $p \nmid u$, because a non-zero $u$ and $\left(x_{ij} - y_{ij}\right)$ both are less than $p$, where $p$ is prime. Therefore, $p \mid \left(x_{ij} - y_{ij}\right) \Rightarrow \left[x_{ij}\,(\text{mod } p)\right] \cong \left[y_{ij}\,(\text{mod } p)\right] \Rightarrow X = Y$ and thus $G_{pAG}(t, u)$ is left cancellative. As every left cancellative AG-groupoid is the right cancellative AG-groupoid [17]. Hence, for prime $p$ and $u \in Z_p - \{0\}$, $G_{pAG}(t, u)$ is a cancellative AG-groupoid.

**Example-10.** $G_{5AG}(3, 4)$ is a cancellative AG-groupoid.

**Theorem-6.** An AG-groupoid of matrices over $Z_n$; $G_{nAG}(t, u)$ is an AG-band, if

$$t + u \cong 1(\text{mod } n).$$

**Proof.** Let $t + u = 1$, to show that a matrix AG-groupoid $G_{nAG}(t, u)$ is an AG-band it is sufficient to show that $A * A = A$;

$$\begin{aligned} A * A &= \left[a_{ij}\right] * \left[a_{ij}\right] \\ &\cong \left[\left(ta_{ij} + ua_{ij}\right)(\text{mod } n)\right] \\ &\cong \left[(t + u)a_{ij}\,(\text{mod } n)\right] \\ &\cong \left[a_{ij}\,(\text{mod } n)\right] \quad \text{as } t + u \cong 1(\text{mod } n) \\ &= A. \end{aligned}$$

Hence a matrix AG-groupoid $G_{nAG}(t, u)$ is an AG-band, if $t + u \cong 1(\text{mod } n)$.

**Example-11.** $G_{5AG}(2, 4)$ is an AG-band.

## 4. EXISTENCE OF AG-GROUP OF MATRICES OVER $Z_n$

In this section, we introduce another class of groupoid of matrices as an AG-groups of matrices over $Z_n$. We study this AG-group of matrices over $Z_n$ and obtain different results. The following theorem shows the existence of AG-groups of matrices over $Z_n$ for $n \geq 3$, and indeed it gives a simple way of construction for matrix AG-groups (mod $n$) of any finite order.

The following theorem guarantees the existence of at least one AG-groups of matrices over $Z_n$ for $n \geq 3$, if $t^2 \cong 1(\text{mod } n)$.

**Theorem-7.** A groupoid of matrices over $Z_n$; $G_n(t, u)$ is an AG-groups of matrices over $Z_n$, if $t^2 \cong 1(\text{mod } n)$ for $t \in Z_n - \{0\}$.

**Proof.** Given that a groupoid of matrices over $Z_n$, $G_n(t, u)$ satisfies $t^2 \cong 1(\text{mod } n)$ for $t \in Z_n - \{0\}$, we have to show $G_n(t, u)$ is matrix AG-group $(\text{mod } n)$.

*Left invertive law:* We show that $(A*B)*C = (C*B)*A$, holds for all $A, B, C \in G_n(t, u)$, since

$$(A*B)*C \cong \left(t^2 a_{ij} + tb_{ij} + c_{ij}\right)(\text{mod } n) \tag{11}$$

and

$$(C*B)*A \cong \left(t^2 c_{ij} + tb_{ij} + a_{ij}\right)(\text{mod } n) \tag{12}$$

This implies that $G_n(t, u)$ is an AG-groupoid $(\text{mod } n)$, as (11, 12) coincide for $t^2 \cong 1(\text{mod } n)$.

*Nonassociativity:* since

$$\text{A} * (\text{B} * \; A*(B*C) \cong \left[\left(ta_{ij} + tb_{ij} + c_{ij}\right)(\text{mod } n)\right]. \tag{13}$$

From (11, 13) it is clear that $G_n(t, u)$ is nonassociative in general.

*Existence of left identity:* '$0 = \left[0_{ij}\right]$' is the left identity of $G_n(t, u)$;

$$0*X \cong \left[x_{ij}(\text{mod } n)\right] = X \; ; \forall \; X \in G_n(t, u)$$

but

$$X*0 \cong \left[\left(tx_{ij}\right)(\text{mod } n)\right] \neq X \text{ in general.}$$

*Existence of inverses:* $(n-1)tX$ or $-tX$ is the inverse of $X \;\forall\; X \in G_n(t, u)$;

$$\begin{aligned}(-tX)*X &= \left[-tx_{ij}\right]*\left[x_{ij}\right] \\ &\cong \left[\left(t\left(-tx_{ij}\right) + x_{ij}\right)(\text{mod } n)\right] \\ &\cong \left(-\left(t^2-1\right)x_{ij}\right)(\text{mod } n) \\ &\cong \left[0_{ij}(\text{mod } n)\right]\end{aligned}$$

and

$$\begin{aligned}X*(-tX) &= \left[x_{ij}\right]*\left[-tx_{ij}\right] \\ &\cong \left[\left(tx_{ij} + \left(-tx_{ij}\right)\right)(\text{mod } n)\right] \\ &\cong \left[0_{ij}(\text{mod } n)\right]\end{aligned}$$

Hence $G_n(t, u)$ is an AG-group of matrices over $Z_n$.

For varying values of $t$ we get different classes of AG-group of matrices over $Z_n$.

**Corollary-3.** Let $G_n(t, u)$ be a groupoid of matrices, then $G_n(n-1, 1)$ is an AG-group of matrices over $Z_n$.

**Proof.** Since $(n-1)^2 \cong 1(\text{mod } n)$. The proof now follows by the Theorem 7.

## 5. CONCLUSION

In this paper a new class of AG-groupoids and AG-groups of modulo matrices over $Z_n$ is investigated, moreover, various types of construction have been introduced for these AG-groupoids and AG-groups. Sufficient examples to show the existence of these notions are provided. It is to be noted that the provided examples generated and verified by various computer programs. The paper contains various nice results, the main result shows that a groupoid is an AG-groupoid of modulo matrices over $Z_n$, if $t^2 \cong u \pmod{n}$ for all $t, u \in Z_n$. Various other relations of these AG-groupoids and AG-groups of modulo matrices over $Z_n$ with some of the already known classes of AG-groupoids are investigated. In Theorem-7, the class of AG-groupoids of modulo matrices over $Z_n$ is further restricted to introduce AG-groups of modulo matrices over $Z_n$. In future, fuzzification of these AG-groupoids and AG-groups of modulo matrices over $Z_n$ will be a nice work.

## ACKNOWLEDGEMENT

The authors are very grateful to Dr. Muhammad Shah, Dr. Muhammad Yousaf and to the unknown referees for their helpful suggestions and comments in order to improve the paper.

No. 44, pp. 2187-2194, 2011.

E-mail address : amanswt@hotmail.com
E-mail address : rashad@uom.edu.pk
E-mail address : iahmaad@hotmail.com
[1] Department of Mathematics University of Malakand, Chakdara Dir L——ower, Pakistan.